\documentclass[12pt,a4paper]{amsart}
\usepackage{geometry}
\usepackage{tabls}
\usepackage{url}
\usepackage[utf8]{inputenc}
\usepackage{amsfonts}
\usepackage{amssymb}
\usepackage{mathrsfs}
\usepackage{amsmath}
\usepackage{amsthm}
\usepackage{amscd}
\usepackage{fancyhdr}
\usepackage{enumerate}
\usepackage{paralist}
\usepackage{xypic}
\usepackage{graphicx}
\usepackage{cite}
\usepackage{lipsum}
\usepackage{footmisc}

\usepackage[all,cmtip]{xy}

\numberwithin{equation}{section}

\theoremstyle{plain}

\newtheorem{Def}[equation]{Definiton}
\newtheorem{Thm}[equation]{Theorem}
\newtheorem{lem}[equation]{Lemma}
\newtheorem{prop}[equation]{Proposition}
\newtheorem{rem}[equation]{Remark}

\newtheorem{conj}[equation]{Conjecture}

\begin{document}

\title{Depth-graded motivic lie algebra}

\author{Jiangtao Li}
\email{ljt-math@pku.edu.cn}
\address{Jiangtao Li \\School of Mathematic Sciences,
        Peking University,
         Beijing, China}

\begin{abstract}

   Consider the neutral Tannakian category mixed Tate motives over $\mathbb{Z}$, in this paper we suggest a way to understand the structure of depth-graded motivic Lie subalgebra generated by the depth one part. We will show that from an isomorphism conjecture proposed by K. Tasaka we can deduce the F. Brown's matrix conjecture and the non-degenerated conjecture about depth-graded motivic Lie subalgebra generated by the depth one part.
\end{abstract}
\maketitle

\section{Introduction}\label{int}

    Denote by  $\zeta_N$  an $N$-th primitive root of unity. Let $\mathcal{MT}(\mathbb{Z}[\zeta_N][1/N])$ denote the category of mixed Tate motives unramified over $\mathbb{Z}[\zeta_N][1/N]$. By the main result of \cite{dg}, the motivic fundamental groupoid of $\mathbb{P}^1-\{0,\mu_N,\infty\}$ can be realized in the category of $\mathcal{MT}(\mathbb{Z}[\zeta_N][1/N])$.

   We call the Lie algebra of the maximal pro-unipotent subgroup of motivic fundamental group of $\mathcal{MT}(\mathbb{Z}[\zeta_N][1/N])$ the motivic Lie algebra of $\mathcal{MT}(\mathbb{Z}[\zeta_N][1/N])$.
    From Proposition 2.3 in \cite{dg}, the motivic Lie algebra  of $\mathcal{MT}(\mathbb{Z}[\zeta_N][1/N])$ is  graded free Lie algebra.

     Since the sub-Tannakian category generated by the function ring of the motivic fundamental groupoid of $\mathbb{P}^1-\{0,\mu_N,\infty\}$ is $\mathcal{MT}(\mathbb{Z}[\zeta_N][1/N]$ for $N=1$ by \cite{brown} and for $N=2,3,4,6,8$ by \cite{del}, from \cite{dg} we know that the motivic Lie algebra of $\mathcal{MT}(\mathbb{Z}[\zeta_N][1/N])$ has an induced depth filtration for $N=1,2,3,4,6,8$.

     In \cite{del}, P. Deligne proves that the depth-graded motivic Lie algebra  of $\mathcal{MT}(\mathbb{Z}[\zeta_N][1/N])$ is a free Lie algebra bi-graded by weight and depth for $N=2,3,4,6,8$. While the structure of depth-graded motivic Lie algebra of $\mathcal{MT}(\mathbb{Z})$ is not fully understood up to now.

     L. Schneps gives the structure of depth-graded motivic Lie algebra of $\mathcal{MT}(\mathbb{Z})$ in depth two \cite{sch}. And A. B. Goncharov's work \cite{gon} give the structure of depth-graded motivic Lie algebra of $\mathcal{MT}(\mathbb{Z})$ in depth three.
      F. Brown gives some conjectural description of the structure of the depth-graded motivic Lie algebra of of $\mathcal{MT}(\mathbb{Z})$ in all depth in \cite{depth}.

      It's widely believed that the Lie subalgebra of depth-graded motivic Lie algebra generated by the depth one part only has the period polynomial relations in depth two among the generators (in \cite{li} we call this statement the non-degenerated conjecture). In this paper, we will show that from an isomorphism conjecture of K.Tasaka \cite{tasaka} we can  deduce Brown's matrix conjecture and the non-degenerated conjecture. Thus we reduce the well-konwn non-degenerated conjecture to a purely linear algebra problem which  probably are more easy to handle.

     And from the analysis in \cite{en}, our results give partial evidence to Brown's homological conjecture about depth-graded motivic Lie algebra in \cite{depth}.
     
\section{Mixed Tate motives}\label{mtm}

Denote $\mathcal{MT}(\mathbb{Z})$ the category of mixed Tate motives over $\mathbb{Z}$. The references about mixed Tate motives are \cite{bf}, \cite{dg}. $\mathcal{MT}(\mathbb{Z})$ is a neutral Tannakian category over $\mathbb{Q}$. Denote $\pi_1(\mathcal{MT}(\mathbb{Z}))$ the fundamental group $\mathcal{MT}(\mathbb{Z})$, then we have
\[
\pi_1(\mathcal{MT}(\mathbb{Z}))=\mathbb{G}_m\ltimes U.
\]
Where $U$ is pro-unipotent algebraic group with free Lie algebra generated by the formal symbol $\sigma_{2n+1}$ in weight $2n+1$ for $n\geq 1$.

By \cite{dg}, the motivic fundamental groupoid of $\mathbb{P}^1-\{0,1,\infty\}$ can be realized in the category $\mathcal{MT}(\mathbb{Z})$.

Denote ${}_0\Pi_1$ the motivic fundamental groupoid of $\mathbb{P}^1-\{0,1,\infty\}$ from the tangential base point $\overrightarrow{1}_0$ at $0$ to the tangential base point $\overrightarrow{-1}_1$ at $1$.  Its function ring over $\mathbb{Q}$ is
\[\mathcal{O}({}_0\Pi_1)=\mathbb{Q}\langle e^0,e^1\rangle,\]
where $\mathbb{Q}\langle e^0,e^1\rangle$ is equipped with the shuffle product.

Denote by ${}_x\Pi_y$ the de-Rham realization of motivic fundamental groupoid of  $\mathbb{P}^1\backslash\{0,1,\infty\}$ from $x$ to $y$ where $x,y\in \{\overrightarrow{1}_0,\overrightarrow{-1}_1\}$. We write $\overrightarrow{1}_0$, $\overrightarrow{-1}_1$ as $0$, $1$ respectively for short. Denote by $G$
 the group of automorphisms of the groupoid ${}_x\Pi_y$ for $x,y\in{0,1}$ which respect to the following structures:

(1) (Groupoid structure) The composition maps
\[
   {}_x\Pi_y \times {}_y\Pi_z\rightarrow {}_x\Pi_z
   \]
for all $x,y,z\in\{0,1\}$.

(2) (Inertia) The automorphism fixes the elements
\[\mathrm{exp}(e_0)\in {}_0\Pi_0(\mathbb{Q}),\, \mathrm{exp}(e_1)\in {}_1\Pi_1(\mathbb{Q}),\]
where $e_0, e_1$ respectively denotes the differential $\frac{dz}{z},\frac{dz}{1-z}$.

From Proposition 5.11 in \cite{dg}, it follows that ${}_x\Pi_y$ is an $G$-torsor.
We have a natural morphism
\[
\varphi:\mathrm{U}^{dR}\rightarrow G\simeq{}_0\Pi_1.
\]

From  \cite{brown} $\varphi$ is injective. Denote by $\mathfrak{g}$ the corresponding Lie algebra of $\mathrm{U}^{dR}$, we have an injective map
\[
i:\mathfrak{g}\rightarrow \mathrm{Lie}\, G \simeq(\mathbb{L} (e_0,e_1),\{ \; , \;\}).
\]
Where ($\mathbb{L}(e_0,e_1),\{\;,\;\}$) is the free Lie algebra generated by $e_0,e_1$ with the following Ihara Lie bracket
\[
\{f,g\}=[f,g]+D_f(g)-D_g(f)
\]
and $D_f$ is a derivation on $\mathbb{L}(e_0,e_1)$ which satisfies $D_f(e_0)=0,\;D_f(e_1)=[e_1,f]$ for $f\in \mathbb{L}(e_0,e_1)$.

We denote by $\mathfrak{h}$ the Lie algebra ($\mathbb{L}(e_0,e_1),\{\;,\;\}$) for short. There is a natural decreasing depth filtration on $\mathfrak{h}$ defined by
\[
\mathfrak{D}^r\mathfrak{h}=\{\xi\in \mathfrak{h}\mid \mathrm{deg}_{e_1}\;\xi\geq r\}.
\]
And define the weight grading by the total degree of $e_0,e_1$ for the elements of $\mathfrak{h}$.
From the injective map $i$, there is an induced depth filtration on $\mathfrak{g}$, define
\[
\mathfrak{dg}=\oplus_{r\geq 1}\mathfrak{D}^r\mathfrak{g}/\mathfrak{D}^{r+1}\mathfrak{g}
\]
with induced Lie bracket as depth graded motivic Lie algebra of $\mathcal{MT}(\mathbb{Z})$.
By Th$\acute{e}$or$\grave{e}$me 6.8(i) in\cite{dg}, we have
$i(\sigma_{2n+1})=(\mathrm{ad}\,e_0)^{2n}(e_1)+$ terms of degree $\geq 2 $ in $e_1$. So $\mathfrak{dg}_1$ is essentially the $\mathbb{Q}$-linear combination of $\overline{\sigma}_{2n+1}=(\mathrm{ad}\,e_0)^{2n}(e_1)$, $n\geq 1$ in $\mathfrak{h}$.

Here we give the definition of restricted even period polynomial:
\begin{Def}\label{polynomial}
For $N\geq 3$, the restricted even period polynomial of weight $N$ is the polynomial $p(x_1,x_2)$ of degree $N-2$ which satisfies\\
(i) $p(x_1,0)=0$, i.e. p is restricted;\\
(ii) $p(\pm x_1,\pm x_2)=p(x_1,x_2)$, i.e. p is even;\\
(iii) $p(x_1,x_2)+p(x_1-x_2,x_1)-p(x_1-x_2,x_2)=0$.\\
Denote $\mathbb{P}_N$ the set of even restricted period polynomials of weight $N$.
\end{Def}

For $\mathbb{Q}$-vector space, denote by $\mathrm{Lie}(V)$ the free Lie algebra generated by the vector space V. Denote by $\mathrm{Lie}_n(V)$ elements of $\mathrm{Lie}(V)$ with exactly $n$ occurrences of the formal Lie bracket $[\;,\;]$.

 For $n\geq2$, define \[\alpha: \mathbb{P} \otimes\underbrace{\mathfrak{d}\mathfrak{g}_1\otimes\cdots\otimes\mathfrak{d}\mathfrak{g}_1}_{n-2}\rightarrow\mathrm{Lie}_n(\mathfrak{d}\mathfrak{g}_1)\]
 by
\[\alpha:\sum p_{r,s}x_1^{r-1}x_2^{s-1}\otimes \overline{\sigma}_{i_1}\otimes\cdots\otimes\overline{\sigma}_{i_{n-2}}\mapsto  \sum p_{r,s}[\cdots[[\overline{\sigma}_r,\overline{\sigma}_s],\overline{\sigma}_{i_1}],\cdots,\overline{\sigma}_{i_{n-2}}]\]
where $[\,,\,]$ is the formal lie bracket. Denote by $\beta:\mathrm{Lie}_n(\mathfrak{d}\mathfrak{g}_1)\rightarrow \mathfrak{d}\mathfrak{g}_n$  the map that replacing the formal Lie bracket by the induced Ihara bracket.

The following conjecture is well-known.
\begin{conj}\label{non}{(non-degenerated conjecture)}
For $n\geq2$, the following sequence
\[\mathbb{P} \otimes\underbrace{\mathfrak{d}\mathfrak{g}_1\otimes\cdots\otimes\mathfrak{d}\mathfrak{g}_1}_{n-2}\xrightarrow{\alpha}\mathrm{Lie}_n(\mathfrak{d}\mathfrak{g}_1)\xrightarrow{\beta} \mathfrak{d}\mathfrak{g}_n\]
is exact.

\end{conj}

\section{Universal enveloping algebra}\label{uea}

Denote by $\mathcal{U}\mathfrak{h}$ the universal enveloping algebra of $\mathfrak{h}$ and denote by $\mathbb{Q}\langle e_0,e_1\rangle$ the non-commutative polynomial ring in symbol $e_0,e_1$. From Proposition $5.9$ in \cite{dg} we know that $\mathcal{U}\mathfrak{h}$ is isomorphic to $\mathbb{Q}\langle e_0,e_1\rangle$ as a vector space. But the new multiplication structure on $\mathbb{Q}\langle e_0,e_1\rangle$ which are transformed from $\mathcal{U}\mathfrak{h}$ are rather subtle. It's not the usual concatenation product .

Denote $gr_{\mathfrak{D}}^r\mathbb{Q}\langle e_0,e_1\rangle$ the elements of $\mathbb{Q}\langle e_0, e_1\rangle$ with exactly $r$ occurrences of $e_1$. We have the following map:
\[
\begin{split}
\rho:&gr_{\mathfrak{D}}^r\mathbb{Q}\langle e_0,e_1\rangle\rightarrow \mathbb{Q}[y_0,y_1,\cdots,y_r]\\
&e_0^{a_0}e_1e_0^{a_1}e_1\cdots e_1e_0^{a_r}\mapsto y_0^{a_0}y_1^{a_1}\cdots y_r^{a_r}
\end{split}
\]
The map $\rho$ is the polynomial representation of $\mathbb{Q}\langle e_0, e_1\rangle$ defined by F. Brown.

In \cite{depth}, F. Brown introduced a $\mathbb{Q}$-bilinear map $\underline{\circ}:\mathbb{Q}\langle e_0, e_1\rangle\otimes_{\mathbb{Q}}\mathbb{Q}\langle e_0, e_1\rangle\rightarrow \mathbb{Q}\langle e_0, e_1\rangle $  which in the polynomial representation can be written as
\[
\begin{split}
f\underline{\circ}g(y_0,\cdots,&y_{r+s})=\sum_{i=0}^sf(y_i,y_{i+1},\cdots,y_{i+r})g(y_0,\cdots,y_i,y_{i+r+1},\cdots, y_{r+s})+\\
&(-1)^{deg f+r}\sum_{i=1}^s f(y_{i+r},\cdots, y_{i+1}, y_i)g(y_0,\cdots,y_{i-1},y_{i+r},\cdots, y_{r+s})
\end{split}
\]
for $f\in \mathbb{Q}[y_0,\cdots,y_r]=\rho(gr_{\mathfrak{D}}^r\mathbb{Q}\langle e_0, e_1\rangle)$, $g\in \mathbb{Q}[y_0,\cdots,y_s]=\rho(gr_{\mathfrak{D}}^s\mathbb{Q}\langle e_0, e_1\rangle)$.

    Since by the general theory of Lie algebra, the natural action of $\mathfrak{h}$ on $\mathcal{U}\mathfrak{h}$ is the form $(a, b_1\otimes b_2\otimes\cdots b_r)\mapsto a\otimes b_1\otimes \cdots \otimes b_r$ in $\mathcal{U}\mathfrak{h}$ for $a, b_1,\cdots, b_r\in \mathfrak{h}$. By Proposition 2.2 in \cite{depth}, we have
 \[
 a_1\circ a_2\circ\cdots\circ a_r=a_1\underline{\circ}(a_2\underline{\circ}(\cdots(a_{r-1}\underline{\circ} a_r)\cdots))
 \]
 for $a_i\in \mathfrak{h}\in \mathbb{Q} \langle e_0, e_1\rangle,1\leq i\leq r-1$, $a_r\in \mathbb{Q}\langle e_0,e_1\rangle$.

    The above formula is still not enough to give a very clear picture of the new multiplication $\circ$ on $\mathbb{Q}\langle e_0, e_1\rangle$. But it's enough for our purpose.

    We first introduce some notation from K.Tasaka \cite{tasaka}.
Denote by
\[
S_{N,r}=\{(n_1,...,n_r)\in\mathbb{Z}^r\mid n_1+...+n_r=N,n_1,...,n_r\geq3:odd\}.
\]
We write $\overrightarrow{m}=(m_1,...,m_r)$ for short, while \[\mathbf{Vect}_{N,r}=\{(a_{n_1,...,n_r})_{\overrightarrow{n}\in S_{N,r}}\mid a_{n_1,...,n_r}\in \mathbb{Q}\}.\]

For a matrix $P=\left(p\dbinom{m_1,...,m_r}{n_1,...,n_r}\right)_{\substack {\overrightarrow{m}\in S_{N,r}\\\overrightarrow{n}\in S_{N,r}}}$, the action of $P$ on $a=(a_{m_1,...,m_r})_{\overrightarrow{m}\in S_{N,r}}$ means
\[
aP=\left(\sum_{\overrightarrow{m}\in S_{N,r}}a_{m_1,...,m_r}p\dbinom{m_1,...,m_r}{n_1,...,n_r}\right)_{\overrightarrow{n}\in S_{N,r}}
\]
Denote by $\mathbb{P}_{N,r}$ the $\mathbb{Q}$-vector space spanned by the set
\[
\{x_1^{n_1-1}\cdots x_r^{n_r-1}\mid(n_1,...,n_r)\in S_{N,r}\}.
\]
Obviously there is an isomorphism
\[\pi:\mathbb{P}_{N,r}\longrightarrow\mathbf{Vect}_{N,r}\]
\[\sum_{\overrightarrow{n}\in S_{N,r}}a_{n_1,...,n_r}x^{n_1-1}_1\cdots x^{n_r-1}_r\longmapsto(a_{n_1,...,n_r})_{\overrightarrow{n}\in S_{N,r}}.\]
Denote by
\[
\mathbf{W}_{N,r}=\{p\in\mathbb{P}_{N,r}\mid p(x_1,...,x_r)=p(x_2-x_1,x_2,x_3,...,x_r)-p(x_2-x_1,x_1,x_3,...,x_r)\}.
\]
Denote by
\[
e\dbinom{m_1,...,m_r}{n_1,...,n_r}=\delta\dbinom{m_1,...,m_r}{n_1,...,n_r}+\sum_{i=1}^{r-1}\delta\dbinom{m_2,...,m_i,m_{i+2},...,m_r}{n_1,...,n_{i-1},n_{i+2},...,n_r}b^{m_1}_{n_i,n_{i+1}}.
\]
Which the $b^m_{n,n'}$ are defined by
\[
b^m_{n,n'}=(-1)^n\dbinom{m-1}{n-1}+(-1)^{n'-m}\dbinom{m-1}{n'-1}.
\]
and $\delta\dbinom{m_1,...,m_r}{n_1,...,n_r}=1$ if $\overrightarrow{m}=\overrightarrow{n}$, $\delta\dbinom{m_1,...,m_r}{n_1,...,n_r}=0$ if $\overrightarrow{m}\neq\overrightarrow{n}$.

And the matrix $E^{(r-i)}_{N,r},\,i=0,1,...,r-2$ are defined by
\[
E^{(r-i)}_{N,r}=\left(\delta\dbinom{m_1,...,m_i}{n_1,...,n_i}e\dbinom{m_{i+1},...,m_r}{n_{i+1},...,n_r}\right)_{\substack {\overrightarrow{m}\in S_{N,r}\\\overrightarrow{n}\in S_{N,r}}}.
\]
We write $E^{(r)}_{N,r}$ as $E_{N,r}$.
And denote by
\[
C_{N,r}=E^{(2)}_{N,r}\cdot E^{(3)}_{N,r}\cdots E^{(r-1)}_{N,r}\cdot E_{N,r}
\]
for $r=2$. And denote by $C_{N,r}$ the one row, one column matrix $1$ for $N>1,\mathrm{odd}$, $r=1$.

By \cite{schneps}, we know $\pi(\mathbf{W}_{N,2})=\mathrm{Ker}\;E_{N,2}$.
For $r\geq 3$, K.Tasaka proved that
\[
(\pi(\mathbf{W}_{N,r})(E_{N,r}-I_{N,r})\subseteq \mathrm{Ker}\; E_{N,r},
\]
which $I_{N,r}$ denotes the identity matrix $$\left(\delta\dbinom{m_1,...,m_r}{n_1,...,n_r}\right)_{\substack {\overrightarrow{m}\in S_{N,r}\\\overrightarrow{n}\in S_{N,r}}}.$$

 Furthermore, K. Tasaka proposed the following conjecture.
 \begin{conj}\label{taisha} {(Tasaka conjecture)}
 The linear map $$\eta:\pi(\mathbf{W}_{N,r})\rightarrow \mathrm{Ker}\;E_{N,r}$$
               $$a_{\overrightarrow{m}}\mapsto (a_{\overrightarrow{m}})(E_{N,r}-I_{N,r})$$
is an isomorphism.
 \end{conj}
 In \cite{tasaka}, K.Tasaka suggests a way to prove the injectivity in the above conjecture. But there is a gap in his proof. We prove the injectivity for $r=3$ in \cite{li}.

 In \cite{depth}, F. Brown proposed the following conjecture
 \begin{conj}\label{fbro}{(Brown's matrix conjecture)} The rank of the matrices $C_{N,r}$ satisfy
 \[
 1+\sum_{N,r>0}\mathrm{rank}\;C_{N,r} x^{N}y^r=\frac{1}{1-\mathbb{O}(x)y+\mathbb{S}(x)y^2}
 \]
 \end{conj}

 Now we can state our main result.
 \begin{Thm}\label{im}
 Tasaka conjecture $\Rightarrow$ Brown's matrix conjecture $\Rightarrow$ non-degenerated conjecture.
 \end{Thm}

 \section{Calculation}
    In this section we will prove Theorem \ref{im}. In fact we will prove a little bit more.

 First we will need the following result about Lie algebra.

 \begin{prop}\label{ideal}
 Let $\mathfrak{L}$ be a  Lie algebra over $\mathbb{Q}$, denote by $\mathcal{U}\mathfrak{L}$ its universal envelope algebra. $\mathfrak{M}$ is a Lie ideal in  $\mathfrak{L}$, denote by $\mathcal{U}\mathfrak{L}(\mathfrak{M})$ the two-sided ideal generated by $\mathfrak{M}$ in $\mathcal{U}\mathfrak{L}$. Then we have
 \[
 \mathfrak{L}\cap(\mathcal{U}\mathfrak{L}(\mathfrak{M}) )=\mathfrak{M}.
 \]

 \end{prop}
 \noindent{\bf Proof}:
 We have the following commutative diagram
 \[
 \xymatrix{
   0  \ar[r] & \mathfrak{M} \ar[d] \ar[r] & \mathfrak{L} \ar[d] \ar[r] & \mathfrak{L}/\mathfrak{M} \ar[d] \ar[r] & 0 \\
   0 \ar[r] & \mathcal{U}\mathfrak{L}(\mathfrak{M}) \ar[r] & \mathcal{U}\mathfrak{L} \ar[r] & \mathcal{U}(\mathfrak{L}/\mathfrak{M}) \ar[r] & 0   }
 \]
 And the first row and second row are short exact sequences. By the Poincar\'{e}-Birkhoff-Witt theorem in Lie algebra, we know that the three vertical maps are all injective. $\mathfrak{L}\cap(\mathcal{U}\mathfrak{L}(\mathfrak{M}) )=\mathfrak{M}$ follows by diagram chasing.    $\hfill\Box$\\

    Recall the injective Lie algebra homomorphism in Section \ref{mtm}
 \[
 i:\mathfrak{g}\rightarrow \mathfrak{h}.
 \]
 Since the Lie algebra $\mathfrak{h}$ is bigraded by weight and depth, the map $i$ induces a natural injective Lie algebra homomorphism
 \[
 \overline{i}:\mathfrak{dg}\rightarrow \mathfrak{h}.
 \]
 The maps $i$ and $\overline{i}$ induce the natural injective algebra homomorphisms on enveloping algebra
 \[
 \mathcal{U}i:\mathcal{U}\mathfrak{g}\rightarrow \mathcal{U}\mathfrak{h}=(\mathbb{Q}\langle  e_0,e_1\rangle,\circ)
 \]
 and
 \[
 \mathcal{U}\overline{i}:\mathcal{U}\mathfrak{dg}\rightarrow \mathcal{U}\mathfrak{h}=(\mathbb{Q}\langle  e_0,e_1\rangle,\circ).
 \]

     As $\mathfrak{g}$ is a free Lie algebra generated by elements $\sigma_{2n+1}$ for $n\geq 1$ in weight $2n+1$. We have $\mathcal{U}\mathfrak{g}=\mathbb{Q}\langle \sigma_3, \sigma_5,\cdots, \sigma_{2n+1},\cdots\rangle$ (the non-commutative polynomial ring generated by the symbol $\sigma_{2n+1}$ for $n\geq 1$) with the usual concatenation product.

    For $r=0$, denote by $L_r$ the rational field $\mathbb{Q}$.  For $r\geq 1$, denote by $L_r$ the $\mathbb{Q}$-linear space generated by elements
\[
\overline{\sigma}_{n_1}\circ \overline{\sigma}_{n_2}\circ\cdots\circ\overline{\sigma}_{n_r}=(\mathrm{ad}\;e_0 )^{n_1-1}e_1\circ (\mathrm{ad}\;e_0 )^{n_2-1}e_1\circ\cdots \circ (\mathrm{ad}\;e_0 )^{n_r-1}e_1
\]
in $\mathcal{U}\mathfrak{h}$ for $n_i\geq 3$, odd, $1\leq i\leq r$.

Define the map $B_r:L_1\otimes_{\mathbb{Q}} L_{r-1}\rightarrow L_r$ by
\[
B_r(\overline{\sigma}_{n_1}\otimes (\overline{\sigma}_{n_2}\circ \cdots\circ\overline{\sigma}_{n_r}))=\overline{\sigma}_{n_1}\circ \overline{\sigma}_{n_2}\circ \cdots\circ\overline{\sigma}_{n_r}.
\]
and  define the map $A_r:\mathbb{P}\otimes_{\mathbb{Q}} L_{r-2}\rightarrow L_1\otimes_{\mathbb{Q}} L_{r-1}$ by
\[
A_r:((\sum_{n_1,n_2\geq 3,\mathrm{odd}} p_{n_1,n_2}x_1^{n_1-1}x_2^{n_2-1})\otimes (\overline{\sigma}_{n_3}\circ\cdots\cdots \circ\overline{\sigma}_{n_r}))=\sum_{n_1,n_2\geq 3,\mathrm{odd}}p_{n_1,n_2}\overline{\sigma}_{n_1}\otimes (\overline{\sigma}_{n_2}\circ \cdots\circ\overline{\sigma}_{n_r} ).
\]
We have the following lemma
\begin{lem}\label{equi}
The non-degenerated conjecture for all $r\geq 2$ is equivalent to that the following sequence is exact
\[
0\rightarrow \mathbb{P}\otimes_{\mathbb{Q}} L_{r-2}\xrightarrow{A_r} L_1\otimes_{\mathbb{Q}} L_{r-1}\xrightarrow{B_r} L_r\rightarrow 0.
\]
for all $r\geq 2$.

\end{lem}
\noindent{\bf Proof}: If
\[\tag{1}
\overline{x}=\sum_{\overrightarrow{n}\in S_{N,r}}a_{n_1,\cdots, n_r}\overline{\sigma}_{n_1}\otimes (\overline{\sigma}_{n_2}\circ\cdots\circ \overline{\sigma}_{n_r})\in \mathrm{Ker}\;B_r,
\]
then by definition we will have
\[\tag{2}
\sum_{\overrightarrow{n}\in S_{N,r}}a_{n_1,\cdots,n_r}\sigma_{n_1}\sigma_{n_2}\cdots \sigma_{n_r}\in \mathfrak{D}^{r+1}\mathcal{U}\mathfrak{g}=\mathfrak{D}^{r+1}\mathbb{Q}\langle \sigma_3,\cdots, \sigma_{2n+1},\cdots\rangle.
\]

    Since $\mathcal{U}\mathfrak{g}$ is a non-commutative polynomial ring, from formula $(1)$ and $(2)$ we have
 \[
 \begin{split}
 \sum_{\overrightarrow{n}\in S_{N,r}}&a_{n_1,\cdots,n_r}\sigma_{n_1}\sigma_{n_2}\cdots \sigma_{n_r}\\
 &=\sum_{i=2}^r\sum_{\overrightarrow{m}\in S_{N,r}}b^i_{m_1,\cdots, m_i, m_{i+1},\cdots,m_r}[[\cdots[\sigma_{m_1},\sigma_{m_2}],\cdots], \sigma_{m_i}]\sigma_{m_{i+1}}\cdots \sigma_{m_r}
\end{split}
 \]
 for some $b^i_{\overrightarrow{m}},\overrightarrow{m}\in S_{N,r}$ and
 \[
 \sum_{(m_1,\cdots,m_i)\in S_{N-m_{i+1}-\cdots-m_r}}b^i_{m_1,\cdots,m_i,m_{i+1},\cdots, m_r}[[\cdots[\sigma_{m_1},\sigma_{m_2}],\cdots],\sigma_{m_i}]\subseteq \mathfrak{D}^{i+1}\mathfrak{g}
 \]

    On one hand, if the non-degenerated conjecture is true for all depth, we will have $\overline{x}\subseteq \mathrm{Im}\;A_r$, i.e. $\mathrm{Im}\; A_r=\mathrm{Ker}\;B_r$.
And since it's obvious that $B_r$ is surjective. While $A$ is injective follows from $\mathrm{Im}\;A_r=\mathrm{Ker}\;B_r$ in depth $r-1$ and the fact that $\mathbb{P}\otimes L_1 \cap L_1\otimes \mathbb{P}=\{0\}$. We deduce that from the non-degenerated conjecture for all depth we will have the short exact sequence for all $r\geq 2$.

    On the other hand, if the sequence is exact for all $r\geq 2$, let
     \[\tag{3}
    \overline{x}=\sum_{\overrightarrow{m}\in S_{N,r}}b_{m_1,\cdots, m_r}\{\{\cdots\{\overline{\sigma}_{m_1},\overline{\sigma}_{m_2}\},\cdots\},\overline{\sigma}_{m_r}\}=0
    \]
in $\mathfrak{dg}_r$, which $\{\,,\,\}$ denotes the induced Ihara Lie bracket on $\mathfrak{dg}$.

Then
\[\tag{4}
x=\sum_{\overrightarrow{m}\in S_{N,r}}b_{m_1,m_2,\cdots,m_r}[[\cdots[\sigma_{m_1},\sigma_{m_2}],\cdots],\sigma_{m_r}]\in gr_{\mathfrak{D}}^{r+1}\mathcal{U}\mathfrak{g},
\]
which $[\,,\,]$ denotes the formal Lie bracket on the non-commutative polynomial ring $\mathcal{U}\mathfrak{g}$. Rewrite $x$ as
\[\tag{5}
x=\sum_{\overrightarrow{m}\in S_{N,r}}a_{m_1,m_2,\cdots,m_r}\sigma_{m_1}\sigma_{m_2}\cdots \sigma_{m_r}.
\]
Denote
\[
x_{\mathfrak{d}}=\sum_{\overrightarrow{m}\in S_{N,r}}a_{m_1,m_2,\cdots,m_r}\overline{\sigma}_{m_1}\otimes (\overline{\sigma}_{m_2}\circ \overline{\sigma}_{m_3}\circ\cdots \circ \overline{\sigma}_{m_r}  )\in L_1\otimes_{\mathbb{Q}} L_{r-1},
\]
then from $(3),(4)$ and $(5)$, we have $x_{\mathfrak{d}}\in \mathrm{Ker}\; B_r$.
Since $\mathrm{Im}\;A_r =\mathrm{Ker}\; B_r$, we have
\[
x_{\mathfrak{d}}=\sum_{\overrightarrow{m}\in S_{N,r}}c_{m_1,m_2,m_3,\cdots,m_r}\overline{\sigma}_{m_1}\otimes (\overline{\sigma}_{m_2}\circ \overline{\sigma}_{m_3}\circ\cdots \circ \overline{\sigma}_{m_r})
\]
and $p=\sum\limits_{m_1,m_2\geq 3,\mathrm{odd}}c_{m_1,m_2,m_3,\cdots,m_r}x_1^{m_1-1}x_2^{m_2-1}\in \mathbb{P}$. Let $\iota: \mathbb{P}\rightarrow \mathfrak{g}$ be the map
\[
\iota:\sum_{r,s\geq 3,\;\mathrm{odd}}p_{r,s}x_1^{r-1}x_2^{s-1}\mapsto \sum_{r,s\geq 3,\;\mathrm{odd}}p_{r,s}[\sigma_r,\sigma_s].
\]
From $\mathrm{Im}\;A_i=\mathrm{Ker}\;B_i$ for $i=2,\cdots,r$, we deduce inductively that
\[
x\in \mathcal{U}\mathfrak{g}(\iota(\mathbb{P}) ),
\]
which $\mathcal{U}\mathfrak{g}(\iota(\mathbb{P}) )$ means the two-sided ideal generated by $\iota({\mathbb{P}})$ in $\mathcal{U}\mathfrak{g}$. By Proposition \ref{ideal}, $x$ belongs to the Lie ideal generated by $\iota(\mathbb{P})$ in $\mathfrak{g}$. So from the short exact sequence for all depth we can deduce the non-degenerated conjecture in all depth.
 $\hfill\Box$\\

    The following lemma reduces the non-degenerated conjecture to a dimension conjecture of $L_r$ in each weight $N$ for all $r$.

    \begin{lem}\label{r}
    Denote by $L_{N,r}$ the weight $N$ part of $L_r$, then
    the formula
    \[
    1+\sum_{N,r>0}\mathrm{dim}_{\mathbb{Q}}\;L_{N,r}x^Ny^r=\frac{1}{1-\mathbb{O}(x)y+\mathbb{S}(x)y^2}
    \]
    is equivalent to that the following sequence is exact
\[
0\rightarrow \mathbb{P}\otimes_{\mathbb{Q}} L_{r-2}\xrightarrow{A_r} L_1\otimes_{\mathbb{Q}} L_{r-1}\xrightarrow{B_r} L_r\rightarrow 0.
\]
for all $r\geq 2$. Which $\mathbb{O}(x)=\frac{x^3}{1-x^2}$, $\mathbb{S}(x)=\frac{x^{12}}{(1-x^4)(1-x^6)}$.
    \end{lem}

\noindent{\bf Proof}:
$'\Rightarrow'$
 It's clear that $B_r$ is surjective and $\mathrm{Im}\;A_r\subseteq \mathrm{Ker}\; B_r$ for all $r\geq 2$.
Since
\[
\sum_{N>0}\mathrm{dim}_{\mathbb{Q}}L_{N,1}x^N=\mathbb{O}(x),
\]
from the dimension formula we have
\[
\tag{6}
\mathrm{dim}_{\mathbb{Q}}\;L_{N,r}x^N-\sum_{N>0}\mathrm{dim}_{\mathbb{Q}}\;L_{N,r-1}\cdot\mathbb{O}(x)+\sum_{N>0}\mathrm{dim}_{\mathbb{Q}}\;L_{N,r-2}\cdot\mathbb{S}(x)=0
\]
for all $r\geq 2$.

It's obvious that $A_2$ is injective, $B_2$ is surjective and $\mathrm{Im}\;A_2\subseteq\mathrm{Ker}\;B_2$. So from formula $(6)$ in $r=2$, we have  $\mathrm{Im}\;A_2=\mathrm{Ker}\;B_2$.

Inductively, we can deduce that  $A_r$ is injective from $\mathrm{Im}\;A_{r-1}=\mathrm{Ker}\;B_{r-1}$ and Goncharov's result $\mathbb{P}\otimes_{\mathbb{Q}}L_1\cap L_1\otimes_{\mathbb{Q}}\mathbb{P}=0$.
Then $\mathrm{Im}\;A_r=\mathrm{Ker}\;B_r$ follows from formula $(6)$ and the fact that $A_r$ is injective, $B_r$ is surjective and $\mathrm{Im}\;A_r\subseteq \mathrm{Ker}\; B_r$.

$'\Leftarrow'$
It's clear that
\[
\sum_{N>0}\mathrm{dim}_{\mathbb{Q}}L_{N,1}x^N=\mathbb{O}(x).
\]
Then from the short exact sequence we have
\[
\mathrm{dim}_{\mathbb{Q}}\;L_{N,r}x^N-\sum_{N>0}\mathrm{dim}_{\mathbb{Q}}\;L_{N,r-1}\cdot\mathbb{O}(x)+\sum_{N>0}\mathrm{dim}_{\mathbb{Q}}\;L_{N,r-2}\cdot\mathbb{S}(x)=0
\]
for all $r\geq 2$. So
 \[
    1+\sum_{N,r>0}\mathrm{dim}_{\mathbb{Q}}\;L_{N,r}x^Ny^r=\frac{1}{1-\mathbb{O}(x)y+\mathbb{S}(x)y^2}.
    \]  $\hfill\Box$\\
\begin{rem}\label{i}
In fact, if we only know that
 \[
    1+\sum_{N,r>0}\mathrm{dim}_{\mathbb{Q}}\;L_{N,r}x^Ny^r\geq \frac{1}{1-\mathbb{O}(x)y+\mathbb{S}(x)y^2},
    \]
which inequality means the coefficient of the term $x^Ny^r$ in the left side is bigger than the corresponding coefficient in the right side for all $N,r>0$, then we can still deduce the short exact sequence exactly the same way as in the proof of Lemma \ref{r}.
\end{rem}

    Now we investigate the polynomial representation of $L_{N,r}$. From the main result of Section \ref{uea}, we have
\[
\rho(\overline{\sigma}_{m_1}\circ\overline{\sigma}_{m_2}\circ\cdots\circ \overline{\sigma}_{m_r})=(y_1-y_0)^{m_1-1}\underline{\circ}((y_1-y_0)^{m_2-1}\underline{\circ}(\cdots((y_1-y_0)^{m_{r-1}-1}\underline{\circ}(y_1-y_0)^{m_r-1})))
\]
for $\overrightarrow{m}=(m_1,m_2,\cdots,m_r)\in S_{N,r}$.

For  $\overrightarrow{n}=(n_1,n_2,\cdots,n_r)\in S_{N,r}$, the coefficient of $y_1^{n_1-1}y_2^{n_2-1}\cdots y_r^{n_r-1}$ in $\rho(\overline{\sigma}_{m_1}\circ\overline{\sigma}_{m_2}\circ\cdots\circ \overline{\sigma}_{m_r}) $ is
\[
c\binom{m_1,m_2,\cdots,m_r}{n_1,n_2,\cdots,n_r}.
\]
Which $c\binom{m_1,m_2,\cdots,m_r}{n_1,n_2,\cdots,n_r}$ is the $(m_1,m_2,\cdots,m_r)$-th row, $(n_1,n_2,\cdots,n_r)$-th column term of the matrix $C_{N,r}$.

 Now we have
\begin{prop}\label{rep}
The following map
\[
\widetilde{\eta}:W_{N,r}\rightarrow \mathbf{Vect}_{N,r}
\]
\[
(a_{m_1,\cdots,m_r})_{\overrightarrow{m}\in S_{N,r}}\mapsto \left(\sum_{\overrightarrow{m}\in S_{N,r}}a_{m_1,\cdots,m_r}\delta\binom{m_1}{n_1}e\binom{m_2,\cdots,m_r}{n_2,\cdots,n_r}                     \right)_{\overrightarrow{n}\in S_{N,r}}
\]
satisfy $\widetilde{\eta}(W_{N,r})\subseteq \mathrm{Ker}\;E_{N,r}$. And furthermore, $\widetilde{\eta}(a)+\eta(a)=0$ for $$a=(a_{m_1,\cdots,m_r})_{\overrightarrow{m}\in S_{N,r}}\in W_{N,r}.$$

\end{prop}
\noindent{\bf Proof}:
Consider the natural action of $\mathfrak{dg}$ on $\mathcal{U}\mathfrak{h}=(\mathbb{Q}\langle e_0,e_1\rangle,\circ )$. By the main results of Section \ref{uea}, we know that
\[
\sum_{({m_2,\cdots,m_r)}\in S_{N-n_1,r-1}}a_{n_1,m_2,\cdots,m_r}e\binom{m_2,\cdots,m_r}{n_2,\cdots,n_r}
\]
is the coefficient of $y_2^{n_2-1}\cdots y_r^{n_r-1}$ in the polynomial representation of
\[
\sum_{(m_2,\cdots,m_r)\in S_{N-n_1,r-1}}a_{n_1,m_2,\cdots,m_r}\overline{\sigma}_{m_2}\circ(e_1e_0^{m_3-1}e_1\cdots e_1e_0^{m_r-1} ).
\]
for $(n_2,\cdots,n_r)\in S_{N-n_1,r-1}$. Furthermore,
\[
\sum_{\overrightarrow{m},\overrightarrow{n}\in S_{N,r}}a_{m_1,m_2,\cdots,m_r}\delta\binom{m_1}{n_1}e\binom{m_2,\cdots,m_r}{n_2,\cdots,n_r}e\binom{n_1,n_2,\cdots,n_r}{k_1,k_2,\cdots,k_r}
\]
is the coefficient of $y_1^{k_1-1}y_2^{k_2-1}\cdots y_r^{k_r-1}$  in the polynomial representation of
\[
\sum_{\overrightarrow{m}\in S_{N,r}}a_{m_1,m_2,\cdots,m_r}\overline{\sigma}_{m_1}\circ \overline{\sigma}_{m_2}\circ(e_1e_0^{m_3-1}e_1\cdots e_1e_0^{m_r-1}).
\]

If $a=(a_{m_1,\cdots,m_r})_{\overrightarrow{m}\in S_{N,r}}\in W_{N,r}$, then
\[
\begin{split}
&\;\;\;\sum_{\overrightarrow{m}\in S_{N,r}}a_{m_1,m_2,\cdots,m_r}\overline{\sigma}_{m_1}\circ \overline{\sigma}_{m_2}\circ(e_1e_0^{m_3-1}e_1\cdots e_1e_0^{m_r-1})\\
&=\frac{1}{2}\sum_{\overrightarrow{m}\in S_{N,r}}a_{m_1,m_2,\cdots,m_r}(\overline{\sigma}_{m_1}\circ \overline{\sigma}_{m_2}-\overline{\sigma}_{m_2}\circ \overline{\sigma}_{m_1})\circ(e_1e_0^{m_3-1}e_1\cdots e_1e_0^{m_r-1})\\
&=    \frac{1}{2}\sum_{\overrightarrow{m}\in S_{N,r}}a_{m_1,m_2,\cdots,m_r}\{\overline{\sigma}_{m_1}, \overline{\sigma}_{m_2}\}\circ(e_1e_0^{m_3-1}e_1\cdots e_1e_0^{m_r-1})                        \\
&=0.
\end{split}
\]
So we have
\[
\sum_{\overrightarrow{m},\overrightarrow{n}\in S_{N,r}}a_{m_1,m_2,\cdots,m_r}\delta\binom{m_1}{n_1}e\binom{m_2,\cdots,m_r}{n_2,\cdots,n_r}e\binom{n_1,n_2,\cdots,n_r}{k_1,k_2,\cdots,k_r}=0,
\]
i.e. $\widetilde{\eta}(W_{N,r})\subseteq \mathrm{Ker}\;E_{N,r}$.

Similarly, in order to prove $\widetilde{\eta}(a)+\eta(a)=0$ for $a=(a_{m_1,\cdots,m_r})_{\overrightarrow{m}\in S_{N,r}}\in W_{N,r}$, it suffices to show that the coefficient of the term $y_1^{n_1-1}y_2^{n_2-1}\cdots y_r^{n_r-1}$ in the polynomial representation of
\[
-\sum_{\overrightarrow{m}\in S_{N,r}}a_{m_1,m_2,\cdots,m_r}[\overline{\sigma}_{m_1}\circ(e_1e_0^{m_2-1}\cdots e_1e_0^{m_r-1})-e_1e_0^{m_1-1}e_1e_0^{m_2-1}\cdots e_1e_0^{m_r-1}]
\]
is equal to the coefficient of the term $y_1^{n_2-1}y_2^{n_3-1}\cdots y_{r-1}^{n_r-1}$ in the polynomial representation of
\[\tag{7}
\sum_{(m_2,m_3,\cdots,m_r)\in S_{N-n_1,r-1}}a_{n_1,m_2,\cdots,m_r}\overline{\sigma}_{m_2}\circ(e_1e_0^{m_3-1}\cdots e_1e_0^{m_r-1})
\]
  for all
$\overrightarrow{n}\in S_{N,r}$.

From Proposition $2.2$ in \cite{depth}, we have
\[
\begin{split}
&\;\;\;\;\;\overline{\sigma}_{m_1}\circ(e_1e_0^{m_2-1}\cdots e_1e_0^{m_r-1})\\
&=\overline{\sigma}_{m_1}e_1e_0^{m_2-1}\cdots e_1e_0^{m_r-1}-e_1\overline{\sigma}_{m_1}e_0^{m_2-1}\cdots e_1e_0^{m_r-1}+e_1e_0^{m_2-1}(\overline{\sigma}_{m_1}\circ(e_1e_0^{m_3-1}\cdots e_1e_0^{m_r-1}))
\end{split}
\]
and the $y_0$ part of the polynomial representation of $(7)$ is
\[
f_0(y_0,y_1,\cdots,y_{r-1})=\sum_{(m_2,\cdots,m_r)\in S_{N-n_1,r-1}}a_{n_1,m_2,\cdots,m_r}((y_1-y_0)^{m_2-1}-y_1^{m_2-1} )y_2^{m_3-1}\cdots y_{r-1}^{m_r-1}.
\]
So  we have
\[
\begin{split}
&\;\;\;\sum_{\overrightarrow{m}\in S_{N,r}}a_{m_1,m_2,\cdots,m_r}[\overline{\sigma}_{m_1}\circ(e_1e_0^{m_2-1}\cdots e_1e_0^{m_r-1})
-e_1e_0^{m_1-1}e_1e_0^{m_2-1}\cdots e_1e_0^{m_r-1}]\\
&=\sum_{\overrightarrow{m}\in S_{N,r}}a_{m_1,m_2,\cdots,m_r}[\overline{\sigma}_{m_1}e_1e_0^{m_2-1}\cdots e_1e_0^{m_r-1}-e_1\overline{\sigma}_{m_1}e_0^{m_2-1}\cdots e_1e_0^{m_r-1}\\
&\;\;\;\;\;\;\;\;\;\;\;+e_1e_0^{m_2-1}(\overline{\sigma}_{m_1}\circ(e_1e_0^{m_3-1}\cdots e_1e_0^{m_r-1}))   -e_1e_0^{m_1-1}e_1e_0^{m_2-1}\cdots e_1e_0^{m_r-1}                                    ]
\end{split}
\]

Since the coefficient of $y_1^{n_1-1}y_2^{n_2-1}\cdots y_r^{n_r-1}$ in the polynomial representation of
\[
\sum_{\overrightarrow{m}\in S_{N,r}}a_{m_1,m_2,\cdots,m_r}e_1e_0^{m_1-1}(\overline{\sigma}_{m_2}\circ(e_1e_0^{m_3-1}\cdots e_1e_0^{m_r-1}))
\]
is equal to $(\eta(a))_{n_1,n_2,\cdots,n_r}$ minus the coefficient of $y_1^{n_1-1}y_2^{n_2-1}\cdots y_r^{n_r-1}$ in $$y_1^{m_1-1}f_0(y_1,y_2,\cdots,y_r).$$
And the polynomial representation of
\[
\begin{split}
&\sum_{\overrightarrow{m}\in S_{N,r}}a_{m_1,m_2,\cdots,m_r}[\overline{\sigma}_{m_1}e_1e_0^{m_2-1}\cdots e_1e_0^{m_r-1}-e_1\overline{\sigma}_{m_1}e_0^{m_2-1}\cdots e_1e_0^{m_r-1}\\
&\;\;\;\;\;\;\;\;\;\;\;\;\;\;\;\;\;\; \;\;\;\;\;\;\;\;\;\;-e_1e_0^{m_1-1}e_1e_0^{m_2-1}\cdots e_1e_0^{m_r-1}                                    ]
\end{split}
\]
is
\[
\begin{split}
&\;\;\;\;f_1(y_0,y_1,\cdots,y_r)\\
&=\sum_{\overrightarrow{m}\in S_{N,r}}a_{m_1,m_2,\cdots,m_r}[(y_1-y_0)^{m_1-1}y_2^{m_2-1}-(y_2-y_1)^{m_1-1}y_2^{m_2-1}\\
&\;\;\;\;\;\;\;\;\;\;\;\;\;\;\;\;\;\;\;\;\;\;\;\;\;\;\;\;\;\;\;\;-y_1^{m_1-1}y_2^{m_2-1}]y_3^{m_3-1}\cdots y_r^{m_r-1}
\end{split}
\]

So to prove $\widetilde{\eta}(a)+\eta(a)=0$ it suffices to prove
\[\tag{8}
\begin{split}
\sum_{\overrightarrow{m}\in S_{N,r}}a_{m_1,m_2,\cdots,m_r}[&e_1e_0^{m_1-1}(\overline{\sigma}_{m_2}\circ(e_1e_0^{m_3-1}\cdots e_1e_0^{m_r-1}))\\
&+e_1e_0^{m_2-1}(\overline{\sigma}_{m_1}\circ(e_1e_0^{m_3-1}\cdots e_1e_0^{m_r-1}))]=0.
\end{split}
\]
and
\[\tag{9}
f_1(y_0,y_1,\cdots,y_r)-y_1^{m_1-1}f_0(y_1,y_2,\cdots,y_r)
\; \mathrm{has}\;\mathrm{trivial}\;y_1^{n_1-1}y_2^{n_2-1}\cdots y_r^{n_r-1}\;\mathrm{term}
\]
for $\overrightarrow{n}\in S_{N,r}$.
The the formula $(8)$ follows from $a_{m_1,m_2,m_3,\cdots,m_r}+a_{m_2,m_1,m_3,\cdots,m_r}=0$.
And the formula $(9)$ follows from that $a=(a_{m_1,\cdots,m_r})_{\overrightarrow{m}\in S_{N,r}}\in W_{N,r}$.$\hfill\Box$\\

\begin{rem}
From Proposition \ref{rep} we obtain Tasaka's result \cite{tasaka} \[\eta(W_{N,r})\subseteq \mathrm{Ker}\;(E_{N,r}-I_{N,r})\]
immediately. See the proof Proposition $5.5$ in \cite{li} for a proof of the fact $$\widetilde{\eta}(a)+\eta(a)=0$$ based on explicit matrix calculation. Also see the proof of Lemma $4.9$ in \cite{dmn} for a proof based on polynomial representation.
\end{rem}

Now we can prove our main results.\\
\noindent{\bf Proof of Theorem \ref{im}}: From linear algebra, we have
\[\tag{10}
\mathrm{Ker}\; C_{N,r}\cong\mathrm{Ker}\;(E^{(2)}_{N,r}E^{(3)}_{N,r}\cdots E^{(r-1)}_{N,r})\oplus \mathrm{Im}\;(E^{(2)}_{N,r}E^{(3)}_{N,r}\cdots E^{(r-1)}_{N,r})\cap \mathrm{Ker}\; E_{N,r}.
\]
By the definition of $C_{N,r}$, view $C_{N,r}$ as linear transformation on the vector space $\mathbf{Vect}_{N,r}$, then we have
\[\tag{11}
\mathrm{Ker}(E^{(2)}_{N,r}E^{(3)}_{N,r}\cdots E^{(r-1)}_{N,r})\cong \bigoplus_{m>1,\mathrm{odd}}\mathrm{Ker}\; C_{N-m,r-1}.
\]
If Conjecture \ref{taisha} (Tasaka conjecture) is true, then from Proposition \ref{rep}  and the fact 
\[
\mathrm{Im}\;(E_{N,r}^{(2)}E_{N,r}^{(3)}\cdots E_{N,r}^{(r-2)})\cap W_{N,r}\cong \bigoplus_{m>0,\mathrm{even}}\mathbb{P}_m\otimes _{\mathbb{Q}}\mathrm{Im}\;C_{N-m.r-2}
\]
we have
\[\tag{12}
\mathrm{Im}\;(E^{(2)}_{N,r}E^{(3)}_{N,r}\cdots E^{(r-1)}_{N,r})\cap \mathrm{Ker}\; E_{N,r}\cong \bigoplus_{m>0,\mathrm{even}}\mathbb{P}_{m}\otimes_{\mathbb{Q}}\mathrm{Im}\; C_{N-m,r-2}.
\]

From formula $(10), (11)$ and $(12)$, we have
\[\tag{13}
\begin{split}
\sum_{N,r>0}\mathrm{dim}_{\mathbb{Q}} \mathrm{Ker}\; C_{N,r} x^N y^r=&\sum_{N,r>0}\mathrm{dim}_{\mathbb{Q}} \mathrm{Ker}\; C_{N,r} x^N y^r\cdot {\mathbb{O}(x)y}\\
&\;\;\;\;\;\;\;\;\;\;\;\;\;\;\;+\mathbb{S}(x)y^2\cdot (1+\sum_{N,r>0}\mathrm{dim}_{\mathbb{Q}} \mathrm{Im}\; C_{N,r} x^Ny^r).
\end{split}
\]

From formula $(13)$, we have
\[\tag{14}
1+\sum_{N,r>0}\mathrm{rank}\;C_{N,r}x^Ny^r=\frac{1}{1-\mathbb{O}(x)y+\mathbb{S}(x)y^2}.
\]

The polynomial representation of element $\overline{\sigma}_{m_1}\circ \overline{\sigma}_{m_2}\circ\cdots \circ \overline{\sigma}_{m_r}$ in $L_{N,r}$ for $\overrightarrow{m}\in S_{N,r}$ is
\[\tag{15}
(y_1-y_0)^{m_1-1}\underline{\circ}((y_1-y_0)^{m_2-1}\underline{\circ}(\cdots\underline{\circ}(y_1-y_0)^{m_r-1})\cdots))
\]
And the coefficient of the term $y_1^{n_1-1}y_2^{n_2-1}\cdots y_r^{n_r-1}$ in the formula $(15)$ is \\
the $(m_1,m_2,\cdots, m_r)$-th row, the $(n_1,n_2,\cdots,n_r)$-th column element of the matrix $C_{N,r}$.
So from formula $(14)$ we have
\[
1+\sum_{N,r>0}\mathrm{dim}_{\mathbb{Q}}\;L_{N,r}x^Ny^r\geq \frac{1}{1-\mathbb{O}(x)y+\mathbb{S}(x)y^2}.
\]
From Remark \ref{i} and Lemma \ref{equi}, we have the non-degenerated conjecture.$\hfill\Box$\\

\begin{rem}\label{dc}
Formula $(12)$ is essentially the Conjecture $4.12$ in \cite{dmn}, in the above proof we actually show that formula $(12)$ is a corollary of Tasaka's isomorphism conjecture. And see \cite{li} for the application of non-degenerated conjecture to motivic multiple zeta values.
\end{rem}

\end{document}